\documentclass[MSNbibl,number,citesort,seceqn,dvips]{arxbj}

\aid{0}
\volume{20}
\issue{2}
\pubyear{2014}
\firstpage{958}
\lastpage{978}
\doi{10.3150/13-BEJ511} %kopijuoti is PTS

\makeatletter
\renewcommand{\l}{\ell}
\newproclaim{condition}{Condition}
\newremark{remark}{Remark}
\newtheorem{Corollary}{Corollary}
\newtheorem{theorem}{Theorem}
\newcommand{\rrvert}{\vert}
\newcommand{\llvert}{\vert}
\def\cal{\mathcal}
\newcommand{\eqref}[1]{(\ref{#1})}
\newcommand{\C}{\cal{C}}
\newcommand{\hbo}{{\hat B}_0}
\newcommand{\hvo}{{\hat V}_0}
\newcommand{\hsi}{\hat{\sigma}}
\newcommand{\hthn}{{\hat{\theta}}_{n}}
\newcommand{\hvphi}{\hat{\varphi}}
\newcommand{\nti}{n\rightarrow\infty}
\newcommand{\raw}{\rightarrow}

\newcommand{\te}{\theta}
\newcommand{\si}{\sigma}

\newcommand{\vphi}{\varphi}

\newcommand{\JAB}{\mathrm{JAB}}
\newcommand{\nppi}{\mathrm{NPPI}}
\newcommand{\NPP}{\mathrm{NPPI}}
\newcommand{\BIAS}{\operatorname{BIAS}}
\newcommand{\VAR}{\operatorname{VAR}}

\newcommand{\np}{\mathrm{NPPI}}
\newcommand{\pt}{\propto}

\newcommand{\MSE}{\operatorname{MSE}}

\newcommand{\cov}{\operatorname{Cov}}
\newcommand{\E}{\mathrm{E}}

\newcommand{\HHJ}{\mathrm{HHJ}}
\newcommand{\BK}{\mathrm{BK}}
\newcommand{\PW}{\mathrm{PW}}
\newcommand{\opt}{\mathrm{opt}}
\newcommand{\infi}{\infty}

\newcommand{\bhoi}{\hat{b}^{\opt,\infi}_{m,\HHJ}}

\newcommand{\lhoi}{\hat{\ell}^{\opt,\infi}_{n,\HHJ}}
\makeatother

\begin{document}
\begin{frontmatter}

\title{Convergence rates of empirical
block length selectors for block bootstrap}
\runtitle{Empirical
block length selectors}

\begin{aug}
%%%% inicialai - be tarpu
\author[1]{\fnms{Daniel J.} \snm{Nordman}\corref{}\thanksref{1}\ead[label=e1]{dnordman@iastate.edu}}
\and
\author[2]{\fnms{Soumendra N.} \snm{Lahiri}\thanksref{2}\ead[label=e2]{snlahiri@stat.tamu.edu}}
\runauthor{D.J. Nordman and S.N. Lahiri} %% auto
\address[1]{Department of Statistics, Iowa State University, Ames, Iowa
50011, USA.\\ \printead{e1}}

\address[2]{Department of Statistics, Texas A\&M University, College
Station, Texas
77843, USA.\\ \printead{e2}}
\end{aug}

% HISTORY:
\received{\smonth{2} \syear{2012}}
\revised{\smonth{9} \syear{2012}}

% ABSTRACT
%
\begin{abstract}
We investigate the accuracy of two general non-parametric methods for
estimating optimal block lengths for block bootstraps with time series --
the first proposed
in the seminal paper of Hall, Horowitz and Jing
(\textit{Biometrika} \textbf{82} (1995) 561--574) and the second from
Lahiri \textit{et al.}
(\textit{Stat. Methodol.} \textbf{4} (2007) 292--321).
The relative performances of these general methods
have been unknown and, to provide a comparison,
we focus on rates of convergence for these
block length selectors
for the moving
block bootstrap (MBB) with variance estimation problems
under the smooth function model.
It is shown that,
with suitable choice of tuning parameters,
the optimal convergence rate of
the first method is $O_p(n^{-1/6})$
%as $n\rightarrow\infty$,
where $n$ denotes the sample size.
The optimal convergence rate of the second method,
with
the same number of tuning parameters,
is shown to be $O_p(n^{-2/7})$, suggesting that the second method
may generally have better large-sample properties for block selection
in block bootstrap applications beyond variance estimation. We also
compare the two general methods
with other plug-in methods specifically designed for block selection in
variance estimation,
where the best possible convergence rate is shown to be $O_p(n^{-1/3})$
and achieved by
a method from Politis and White
(\textit{Econometric Rev.} \textbf{23} (2004) 53--70).
%A key result is a uniform approximation for the
%mean square error function of the MBB variance estimators
%over a range of block sizes that may be of independent
%interest.
\end{abstract}

% KEYWORDS
% visi is mazosios raides ir pagal abecele
%

\end{frontmatter}

%s1 #&#
\section{Introduction}
\label{sec1}
Performance of block bootstrap methods critically depends on
the choice of block lengths. A~common approach to the problem
is to choose a block length that minimizes the Mean Squared Error
(MSE) function of block bootstrap estimators as a function of
the block length. For many important functionals, expansions
for the MSE-optimal
block lengths are known.
If $\hat{\theta}_n$ denotes an estimator of a parameter of interest
$\theta\in\mathbb{R}$
based on a stationary stretch $X_1,\ldots,X_n$, examples of relevant
functionals $\varphi_n$
of the distribution of $\hat{\theta}_n$ include the bias $\varphi_{1n}=
\E(\hat{\theta}_n-\theta)$,
variance $\varphi_{2n}= \operatorname{Var}(\hat{\theta}_n)$, and
the distribution
function $\varphi_{3n}(x_0)=
P( \sqrt{n}(\hat{\theta}_n-\theta)/\tau_n \leq x_0 )$ (i.e., given
$x_0\in\mathbb{R}$ and where $\tau_n^2$
represents either the variance of $\sqrt{n}(\hat{\theta}_n-\theta)$ or
an estimator of this, cf.~\cite{L07}).
If $\hat{\varphi}_n(\ell)$ denotes a block bootstrap estimator of
$\varphi_n$ based on block length $\ell$, then
as $n\rightarrow\infty$ the bias and variance of $\hat{\varphi
}_n(\ell
)$ often admit
expansions of the form
%e1.1 #&#
%
\begin{equation}
\label{eqngenexp} n^{2a} \operatorname{Var}\bigl( \hat{\varphi}_n(
\ell) \bigr) =V_0 \frac{\ell^r}{n} \bigl(1+o(1) \bigr),\qquad
n^a \operatorname{Bias} \bigl( \hat{\varphi}_n(\ell)
\bigr) = -\frac{B_0}{\ell} \bigl(1+o(1) \bigr)
\end{equation}
for some known constants $a,r>0$ depending on $\varphi_n$ (e.g.,
$a=r=1$ for functionals
$\varphi_n= \varphi_{1n},\varphi_{2n}$, while $r=2,a=1/2$ for the
distribution function $\varphi_n= \varphi_{3n}(x_0)$ when $|x_0|\neq1$)
and lead to a large sample approximation of MSE-optimal block size
given by
%e1.2 #&#
%
\begin{equation}
\label{eqngenblock} \ell_n^0\equiv
\ell_n^0(\varphi) = \mathcal{C}_0 n
^{{1}/{(r+2)}} \bigl(1+o(1) \bigr),\qquad \mathcal{C}_0 \equiv\biggl(
\frac{2 B_0^2}{ r V_0} \biggr)^{{1}/{(r+2)}},
\end{equation}
involving population quantities $B_0=B_0(\varphi_n), V_0=V_0(\varphi_n)
\in\mathbb{R}$
that depend on the functional~$\varphi_n$, the bootstrap method, and
various parameters
of the underlying process. For smooth function model statistics $\hat
{\theta}_n$ (described below),
these expansions (\ref{eqngenexp}) have been established
for the moving block and non-overlapping block methods \cite
{H95,K89,L99,L07} and, in
particular, are also known for the variance functional $\varphi_{2n}$
with other block
bootstraps, such as the circular block bootstrap \cite{PR92} and
stationary bootstrap \cite{N09,PR94};
see \cite{L03} and references therein.
However, as the theoretical approximations (\ref{eqngenblock}) for
the optimal
block lengths typically depend on different unknown population
parameters of the underlying process in an intricate manner, these are
not directly usable in practice.

Different data-based methods for the selection of optimal block
lengths have been proposed in the literature. One of the most popular
general methods is proposed by Hall, Horowitz and Jing~\cite{H95}
(hereafter referred to as HHJ) which employs a
subsampling method (cf. \cite{PR94}) to construct
an empirical version of the MSE function and minimizes this
to produce an estimator of the optimal block length.
We will refer to this approach as the HHJ method.
A second general method for selecting the optimal
block length is put forward by Lahiri \textit{et al.} \cite{L07}.
This method is based on the jackknife-after-bootstrap method
of Efron \cite{E92} and its extension to block bootstrap by
Lahiri \cite{L02}. For reasons explained in \cite{L07}
(see also Section \ref{sec2} below), we will refer to this
method as the non-parametric plug-in method (or the NPPI method,
in short). Both the HHJ and NPPI methods
are called ``general'' because these can be used in the same manner
across different functionals (e.g., bias, variance,
distribution function, quantiles, etc.) to find the optimal
block size for bootstrap estimation, \textit{without} requiring exact
analytical expressions for the corresponding optimal block length approximation
(\ref{eqngenblock}) (i.e., without requiring explicit forms for
quantities $B_0,V_0$).
In particular,
for a given functional,
the HHJ method aims to directly estimate the constant $\mathcal{C}_0$
in the optimal
block approximation (\ref{eqngenblock}) while the NPPI
method separately and non-parametrically estimates the bias $B_0$ and
variance $V_0$
quantities in (\ref{eqngenblock}) without structural knowledge of these.
Our major objective here is to investigate the convergence rates
of these two \textit{general} methods.
For instance, despite the popularity of the HHJ method,
little is theoretically known about its properties
for block selection or how this compares to the NPPI method.
As a context to compare the methods, we focus on their performance for
block selections
in \textit{variance} estimation problems with the block bootstrap.
In the literature,
a few other block length selection methods also exist. These
are primarily plug-in estimators which \textit{necessarily}
require an explicit expression for the optimal block approximation
(\ref
{eqngenblock})
for each specific functional and
for each block bootstrap method (i.e., requiring exact forms for
$B_0,V_0$) and are not
the focus of this paper. However,
two popular plug-in methods for
the variance functional in the latter category are
given by B\"uhlmann and K\"unsch \cite{BK99} and
Politis and White \cite{PW04} (and its corrected version
Patton, Politis and White \cite{P10}). For completeness,
we later compare the performance of the two general
methods with these plug-in methods
for block selection in variance estimation.
% of B\"uhlmann and K``unsch (1989) and of
% Politis and White (2003) and Politis et al. (2010).

For concreteness, we shall restrict attention to the moving block
bootstrap (MBB) method \cite{K89,L92},
which was the original focus of the HHJ method \cite{H95} and the
plug-in method of B\"uhlmann and K\"unsch \cite{BK99}
and
shares close large-sample connections to other
block bootstrap methods (e.g., circular block bootstrap,
non-overlapping block bootstrap, untapered version of the tapered block
bootstrap)
\cite{L99,N09,PP01,PW04}. Further, we shall
work under the
\textit{smooth function model} of Hall \cite{H92} (see Section \ref
{sec21} below)
which provides a convenient theoretical framework but, at the
same time, is general enough to cover many commonly used estimators
in the time series context (\cite{L03}; Chapter 4).
Accordingly, let $\hthn$ be an estimator
of a parameter of interest $\te$ under the smooth function model
and suppose that the MBB is used for
estimating $\si^2_n \equiv n \operatorname{Var}(\hthn)$
or its limiting form
%e1.3 #&#
%
\begin{equation}
\label{sigma} \sigma_\infty^2 \equiv\lim
_{n \to\infty}n \operatorname{Var}(\hat{\theta}_n).
\end{equation}
Let
%e1.4 #&#
%
\begin{equation}
\label{mse0} \MSE_n(\ell) \equiv\E\bigl\{\hat{
\sigma}^2_{n}( \ell)-\sigma_\infty^2
\bigr\}^2
\end{equation}
denote the MSE of the MBB variance estimator $\hat{\si}_n^2(\ell)$
based on
blocks of length $\ell$ and a sample of size $n$. (Defining the MSE with
$\sigma_n^2$ or $\sigma_\infty^2$ makes no difference in the
following and,
for clarity, it is helpful to fix a target $\sigma_\infty^2$ in
defining (\ref{mse0}) throughout.)

The theoretical MSE-optimal block size
is given by
%e1.5 #&#
%
\begin{equation}
\label{opt} \ell_n^{\opt} = \operatorname{argmin} \bigl\{
\MSE_n(\ell)\dvt\ell\in\mathcal{J}_n \bigr\},
\end{equation}
where $\mathcal{J}_n$ is a suitable set of block lengths including the
optimal block length. As alluded to above (\ref{eqngenexp}),
under some standard regularity conditions, it can be shown that
\[
\MSE_n(\ell) \approx f_n(\ell) \equiv B_0^2\ell^{-2} + V_0
n^{-1}\ell,\qquad
\ell\in\mathcal{J}_n,
\]
where $B_0$ and $V_0$ are population parameters arising, respectively,
from the bias and variance
of the MBB variance estimator $\hat{\si}_n^2(\ell)$.
Let $\ell_n^0 \equiv\operatorname{argmin} \{f_n(\ell) \dvt\ell>0\}
=\mathcal{C}_0 n^{1/3}$
denote the minimizer of the asymptotic approximation $f_n(\cdot)$ to the
MSE function, where $\mathcal{C}_0 = [2 B_0^2/V_0]^{1/3}$ (cf.~(\ref
{eqngenblock})).
As a first step towards investigating the accuracy of
different empirical block rule selection methods,
we consider the relative error of this
theoretical approximation and
show that
\[
\frac{\ell_n^{\opt} -\ell_n^0}{\ell_n^0} = O\bigl(n^{-1/3}\bigr)
\]
as $n\rightarrow\infty$. Thus,
the true optimal block size and the optimal block size
determined by the asymptotic approximation to the MSE curve of
the block bootstrap estimator differ by a margin of
$O(n^{-1/3})$ on the relative scale.
%Further, it is observed that
In general, this rate cannot be improved
further. As a result, for empirical block length selection rules
involving estimation steps
that target $\ell_n^0$ (which all existing methods do), the upper bound
on their accuracy for estimating the true optimal block length
$\ell_n^{\opt}$ is $O_p(n^{-1/3})$.
%On the other hand, the block
%length selection rules that attempt to estimate $\ell_n^{\opt}$
%directly have a slower rate of convergence than $O_p(n^{-1/3})$.
% Consequently, %as an implication of the above result, it appears
%%that
%$O_p(n^{-1/3})$ is the \textit{optimal} achievable rate of
%converegence for a block length selection rule.

Next, we consider the convergence rates of the two general methods.
Let $\hat{\ell}_{n,\HHJ}^{\opt}$
and $\hat{\ell}_{n,\nppi}^{\opt}$, respectively,
denote the estimators of the
optimal block length based on the HHJ
and NPPI methods.
We show that
under some mild conditions and with a suitable choice of the
tuning parameters,
\[
\frac{\hat{\ell}_{n,\HHJ}^{\opt} - \ell_n^{\opt}}{\ell_n^{\opt
}} = O_p\bigl(n^{-1/6}\bigr)
\]
as $n\rightarrow\infty$. Thus, the (relative) rate of convergence of
the HHJ estimator of the optimal block length is $O_p(n^{-1/6})$.
The block length in block bootstrap methodology plays a role similar
to a smoothing parameter in non-parametric functional estimation. It is
well known (cf. \cite{H87}) that non-parametric data based
rules for bandwidth estimation often have an ``excruciatingly slow''
(relative) rate of convergence (e.g., of the order of $O_p(n^{-1/10})$).
The convergence rate of the HHJ method turns out to be
relatively better. It is worth noting that
the HHJ block estimator, based on the overlapping version of
the subsampling method, has the \textit{same} rate of convergence
irrespective of the dependence structure of the underlying
time series $\{X_t\}$. Additionally, in the process of determining this
convergence rate,
we also provide the theoretical guidance on optimally choosing two tuning
parameters required in implementing the HHJ method,
which has been an unresolved aspect of the method.

Next, we consider the NPPI method and compare its relative
performance with the HHJ method. The rate of convergence of the
NPPI method is determined by two factors, which
arise from estimating the variance and the bias
of a block bootstrap estimator (i.e., quantities $V_0$ and $B_0$
appearing in
$\ell_n^0=\mathcal{C}_0 n^{1/3}$, $\mathcal{C}_0 = [2
B_0^2/V_0]^{1/3}$). The factor due to
the variance part is based on the
(block) jackknife-after-bootstrap method \cite{E92,L02}, and it
attains an optimal rate
of $O_p(n^{-2/7})$, with a suitable choice of the
tuning parameters. On the other hand, the second factor
is determined by a non-standard bias estimator that
turns out to be adaptive to the strength of dependence of
$\{X_t\}$. Let $r(k)$ denote the autocovariance function of
(a suitable linear function of) the $X_t$'s. When $r(k) \sim Ck^{-a}$
as $k\raw\infty$ for a suitably large $a>1$,
the rate of convergence of the second term
can be as small as $O_p(n^{-1/2 +\varepsilon})$,
for a given $\varepsilon>0$, with a suitable choice of
the tuning parameters. Thus, combining the
two, the optimal
rate of convergence of the NPPI method
%, based on both factors
becomes $O_p(n^{-2/7})$, which is better than optimal rate
$O_p(n^{-1/6})$ for the HHJ method. For this to hold,
the user needs to specify \textit{two} tuning parameters, the same
number as with the HHJ method.
Also, the convergence rate $O_p(n^{-2/7})$ is interesting in
the variance estimation problem
because this matches the best rate obtained by the plug-in block
selection method of
B\"uhlmann and K\"unsch~\cite{BK99}. Their method is
a four-step algorithm which uses lag weight estimators of the
spectral density at zero and again requires explicit forms for
quantities appearing in the bias and variance (e.g., $B_0,V_0$)
of the MBB variance estimator. Hence, while the NPPI method
for block selection applies more generally to other functionals,
its convergent rate matches the optimal one for a plug-in method
specifically tailored to the variance estimation problem.
This provides some evidence supporting the use of the NPPI method
in block selection with other functionals outside of variance estimation.

The rest of the paper is organized as follows. In Section \ref{sec2}, we
briefly describe the smooth function model, the MBB and the
empirical block length selectors proposed by HHJ \cite{H95}
and Lahiri \textit{et al.} \cite{L07}. In Section \ref{sec3}, we present
the conditions and derive a
general result on uniform approximation of the MSE
of a block bootstrap estimator which may be of
independent interest. We describe main results on the HHJ
and the NPPI methods in Sections \ref{sec4} and \ref{sec5}, respectively.
In Section \ref{sec55}, we compare the general HJJ/NPPI methods
with other plug-in
block selection approaches for the MBB in the variance estimation problem.
In particular, a plug-in method of Politis and White \cite{PW04} (see
also \cite{P10})
is shown to achieve the best possible convergence rate for block selection
with variance functionals.
Section \ref{sec6} sketches proofs of the main results, where full
proofs are
deferred to a supplementary material appendix \cite{NL12}.

%s2 #&#
\section{Preliminaries}
\label{sec2}
%s2.1 #&#
\subsection{MBB variance estimator and optimal block length}
\label{sec21}
Let $\mathcal{X}_n =(X_1,\ldots,X_n) $ be a stationary stretch of
$\mathbb{R}^d$-valued random vectors with mean $\E
X_t=\mu\in\mathbb{R}^d$. We shall consider the problem of
estimating the variance of a statistic framed in the ``smooth
function'' model \cite{H92}. Using some function $H\dvtx\mathbb{R}^d
\rightarrow\mathbb{R}$ and the sample mean $\bar{X}_{n}
=\sum_{i=1}^n X_i/n$, suppose that a statistic can be expressed as
$\hat{\theta}_n= H(\bar{X}_{n})$ for purposes of estimating a
process parameter $\theta= H(\mu)$. The ``smooth
function'' model covers a wide range of parameters and their estimators,
including sample mean, sample autocovariances, Yule--Walker estimators,
among others; see Chapter 4, \cite{L03} for more examples.
Recall the target variance of interest
is $\sigma^2_n\equiv n \operatorname{Var}(\hat{\theta}_n)$ or its
limit (\ref{sigma}).

We next describe the MBB variance estimator. Let $\ell<n\in\mathbb{N}$
(set of positive
integers) denote the block length and create overlapping length
$\ell$ blocks from $\mathcal{X}_n$ as $\{ \mathcal{X}_{i,\ell}\dvt
i=1,\ldots, n-\ell+1\}$, where $\mathcal{X}_{i,\ell}=
(X_i,\ldots,X_{i+\ell-1})$ for any integer $i,\ell\geq1$. We
independently resample $\lfloor n/\ell\rfloor$ blocks by letting
$I_1,\ldots,I_{\lfloor n/\ell\rfloor}$ denote i.i.d. random variables
with a uniform distribution over block indices
$\{1,\ldots,n-\ell+1\}$ and then define a MBB sample $X_1^*,\ldots,
X^*_{n_1}$ of size $n_1=\ell\lfloor n/\ell\rfloor$ as
$(\mathcal{X}_{I_1,\ell},\ldots,
\mathcal{X}_{I_{ \lfloor n/\ell\rfloor},\ell})$, where $\lfloor x
\rfloor$
denotes the integer part of a real number $x$.
The MBB analog of $\hat{\theta}_n$ is given by
$\hat{\theta}_n^* = H(\bar{X}_n^*)$ using the MBB sample mean
$\bar{X}_n^* = \sum_{i=1}^{n_1}X^*_i/n_1$ and the MBB variance
estimator is then defined as
\[
\hat{\sigma}_{n}^2(\ell) \equiv n_1\operatorname{Var}_*
\bigl(\hat{\theta}_n^* \bigr),
\]
where $\operatorname{Var}_*(\cdot)$ denotes the variance with respect
to the
bootstrap distribution conditional on the data $\mathcal{X}_n$.

For variance estimation, we briefly consolidate notation from Section
\ref{sec1} on optimal block lengths.
The
performance of the MBB again depends on the block choice $\ell$.
% and
%the optimal block size $\ell_n^{0}$ for minimizing the
%mean-squared error
% balances a well-known trade-off between the bias and
%variance of the MBB estimator.
Under certain dependence conditions and block
assumptions ($\ell^{-1} + \ell/n\rightarrow0$), the
asymptotic bias and variance of the MBB estimator are
%e2.1 #&#
%
\begin{equation}
\E\hat{\sigma}_{n}^2(\ell) -\sigma_\infty^2
= -\frac{B_0}{\ell} \bigl(1+o(1) \bigr), \qquad\operatorname{Var}\bigl
[\hat{\sigma}_{n}^2(
\ell) \bigr] = V_0\frac{\ell}{n} \bigl(1+o(1) \bigr)\label{bv}
\end{equation}
as $n\rightarrow\infty$, for some population parameters $B_0, V_0$
depending
on the covariance structure
of the underlying process (cf. \cite{H95,K89} and Condition \ref{coS} of
Section \ref{sec31}).
%>From the
%variance and squared bias expressions,
Thus, the main component in MSE~(\ref{mse0}) of the MBB follows as
%e2.2 #&#
%
\begin{equation}
\label{mse11} \MSE_n(\ell) \approx f_n(\ell) \equiv
\frac{B^2_0}{\ell^2} + V_0\frac{\ell}{n}
\end{equation}
as
$n\rightarrow\infty$. The minimizer of $f_n(\ell)$ is given by
%e2.3 #&#
%
\begin{equation}
\l_n^0\equiv\mathcal{C}_0 n^{1/3},
\label{l-0}
\end{equation}
where $\mathcal{C}_0=[2B_0^2/V_0 ]^{1/3}$. From \eqref{mse11} and
\eqref{l-0},
the optimal block minimizing $\MSE_n(\ell)$ behaves as
$\ell_n^{\opt}\approx\l_n^0=\mathcal{C}_0
n^{1/3}$ in large samples \cite{H95,K89,L03}.
As a result, to examine properties of the block length selection
methods, we shall create a
collection of block lengths $\mathcal{J}_n \equiv\{ \ell\in\mathbb
{N} \dvt K^{-1} n^{1/3} \leq\ell\leq K n^{1/3} \}$, for a suitably\vspace*{1pt}
large constant $K>0$ such that $K^{-1}<\mathcal{C}_0 <K$,
and formally define
the optimal block size $\ell_n^{\opt}$ as in (\ref{opt}).

% Note that for any sequence $\ell_n \in\mathcal{J}_n$ (i.e.,
%of order $n^{1/3}$), the error $\MSE_n(\ell_n)$ will have the optimal
%order
%$O(n^{-2/3})$ under appropriate mixing/moment conditions (cf. Theorem?)

%s2.2 #&#
\subsection{The Hall--Horowitz--Jing (HHJ) block estimation method}
\label{sec22}
The HHJ \cite{H95} method seeks to estimate the optimal block size
$\ell_n^{\opt}$ by minimizing an empirical version of
the MSE (\ref{mse0}) created by subsampling (data blocking). Let
$m\equiv m_n \in\mathbb{N}$ denote a
sequence satisfying $m^{-1} + m/n \rightarrow0$
as $n \rightarrow\infty$, which serves to define the length of
subsamples $\mathcal{X}_{i,m}=(X_i,\ldots,X_{i+m-1})$,
$i=1,\ldots,n-m+1$. For each subsample, let
$\hat{\sigma}^2_{i,m}(b)$ denote the MBB variance estimator
resulting from resampling length $b$ blocks from observations
$\mathcal{X}_{i,m}$. For clarity, note that MBB block lengths on
size $m$ subsamples are denoted by ``$b$,'' while ``$\ell$''
denotes MBB block lengths applied to the original data
$\mathcal{X}_n$. To approximate the error $\MSE_{m}(b)\equiv
\E\{\hat{\sigma}^2_{m}(b)-\sigma_\infty^2 \}^2$ in MBB variance
estimation incurred by using length $b$ blocks in samples of size
$m$, we form a subsampling estimator
%e2.4 #&#
%
\begin{equation}
\label{mse} \widehat{\MSE}_m(b) = \frac{1}{n-m+1}\sum
_{i=1}^{n-m+1} \bigl[ \hat{\sigma}^2_{i,m}(b)
- \hat{\sigma}^2_{n}(\tilde{\ell}_n)
\bigr]^2,
\end{equation}
where the initializing MBB estimator
$\hat{\sigma}^2_{n}(\tilde{\ell}_n)$ of $\sigma_\infty^2$ is
based on
the entire sample
$\mathcal{X}_n$ and on
a plausible pilot block size $\tilde{\ell}_n$. By
minimizing $\widehat{\MSE}_m(b)$ over $\mathcal{J}_m$, we formulate
%e2.5 #&#
%
\begin{equation}
\label{b1} \hat{b}^{\opt}_{m,\HHJ}= \operatorname{argmin} \bigl\{
\widehat{\MSE}_m(b)\dvt b \in\mathcal{J}_m \bigr\}
\end{equation}
as an estimator of the theoretically optimal MBB block length
$b_m^{\opt}$ for a size $m$ sample, with
%e2.6 #&#
%
\begin{equation}
\label{b3} b^{\opt}_m = \operatorname{argmin} \bigl\{
\MSE_m(b)\dvt b \in\mathcal{J}_m \bigr\}.
\end{equation}
Next, is a rescaling step that involves approximating true optimal
block length $\ell_n^{\opt}$
with the minimizer $\l_n^0$ of MSE-approximation \eqref{mse11}. That
is, as $b^{\opt}_m$ is the ``size $m$
sample version'' of $\ell_n^{\opt}$ in (\ref{opt}),
one uses the large-sample block approximation $b_m^{\opt}\approx b_m^0
=\mathcal{C}_0
m^{1/3}$ and $\ell_n^{\opt}\approx\l_n^0 =\mathcal{C}_0 n^{1/3}$ from
(\ref{l-0}) to
re-scale $\hat{b}_{m,\HHJ}^{\opt}$ and subsequently define the HHJ
estimator of $\ell_n^{\opt}$ as
%e2.7 #&#
%
\begin{equation}
\label{HHJblock} \hat{\ell}_{n,\HHJ}^{\opt} = (n/m)^{1/3}
\hat{b}_{m,\HHJ}^{\opt}.
\end{equation}

Hence, the HHJ method requires specifying both a subsample size $m$ and
a pilot
MBB block size $\tilde{\ell}_n$, which impact the
performance of the block estimator $\hat{\ell}_{n,\HHJ}^{\opt}$.
%We later provide optimal sizes of these tuning parameters and
%establish an optimal rate of convergence for
%$(\hat{\ell}_{n,\HHJ}^{\opt} -\ell_n^{\opt})/\ell_n^{\opt}$.

%s2.2.1 #&#
\subsubsection{An oracle-like subsampling MSE}
\label{sec221}
For purposes of comparison with the HHJ method,
we also define a second subsampling MSE given as
%e2.8 #&#
%
\begin{equation}
\label{mse2} \widehat{\MSE}_m^{\infi}(b) = \frac{1}{n-m+1}
\sum_{i=1}^{n-m+1} \bigl[ \hat{
\sigma}^2_{i,m}(b) - \sigma^2_\infty
\bigr]^2,
\end{equation}
which resembles the empirical MSE (\ref{mse}) after
replacing the variance estimator $\hat{\sigma}_n^{2}(\tilde{\ell
}_n)$ with
its target $\sigma_\infty^2$ from (\ref{sigma}). This subsampling
MSE serves to remove one tuning parameter $\tilde{\ell}_n$ in the
original HHJ method by unrealistically assuming
$\sigma_\infty^2$ is known. However, we may parallel the performance of
the HHJ block estimators $\hat{b}^{\opt}_{m,\HHJ}$ and $\hat{\ell
}^{\opt
}_{n,\HHJ}$ to their
oracle-like counterparts
%e2.9 #&#
%
\begin{equation}
\label{b2} \bhoi= \operatorname{argmin} \bigl\{ \widehat{\MSE
}_m^{\infi}(b)
\dvt b \in\mathcal{J}_m \bigr\}
\end{equation}
based on (\ref{mse2}) and the resulting estimator of the optimal
block length $\ell_n^{\opt}$ given by%e2.10 #&#
\begin{equation}
\label{HHJblock2} \lhoi= (n/m)^{1/3} \bhoi.
\end{equation}
Both $\hat{\ell}^{\opt}_{n,\HHJ}$ and
$\lhoi$ estimate the same optimal
block size $\ell^{\opt}_n$, but the estimator $\lhoi$
is based on an \textit{unbiased} subsampling criterion
through knowledge of
$\sigma_\infty^2$, that is,
$\E[\widehat{\MSE}_m^{\infi} (b)] =\MSE_m (b)$ for all $b\in
\mathcal{J}_m$.

%s2.3 #&#
\subsection{The non-parametric plug-in (NPPI) method}
\label{sec23}
The NPPI method is based on
the non-parametric plug-in principle \cite{L07}
which yields estimators of MSE optimal smoothing parameters
in general non-parametric function estimation problems.
Here we describe the method
for estimating the optimal block length
for the variance functional using the MBB.
Like any plug-in method, the target quantity for
the NPPI method is the minimizer $\l_n^0$ of
the MSE-approximation
$ f_n(\ell)$ of \eqref{mse11},
which again is of the form
$
\l_n^0= \C_0n^{1/3}
$
from~\eqref{l-0} with population
parameters $B_0$ and $V_0$ in $\C_0= [2 B_0^2/V_0]^{1/3}$ determined
by the
bias and variance expansion \eqref{bv} of the MBB variance estimator.
The NPPI method
estimates the bias and the variance of the
MBB estimator non-parametrically, and then
estimates $B_0$ and $V_0$ by
inverting \eqref{bv}.
Specifically, the method constructs estimators $\widehat{\operatorname
{BIAS}} $
and $\widehat{\operatorname{VAR}}$ satisfying
\[
\frac{\widehat{\VAR}}{\operatorname{Var}(\hat{\si}^2_n(\ell_1))}
\stackrel{p}
{\rightarrow} 1,\qquad \frac{\widehat{\BIAS}}{\operatorname{Bias}(\hat
{\si}^2_n(\ell_2))}
\stackrel{p}
{\rightarrow} 1 \qquad\mbox{as $\nti$}
\]
for some block lengths $\ell_1$ and $\ell_2$
and defines
$
\hat{V}_0= [n\ell_1^{-1}] \widehat{\VAR}
$ and
$
\hat{B}_0= \ell_2 \widehat{\BIAS}
$.
Then, the ${\NPP}$ estimator of the optimal block
length %$\ell^0$
is given by
%e2.11 #&#
%
\begin{equation}
%eqn(2.7)
\hat{\ell}^0_{\np} = \bigl[2
\hat{B}_0^2/\hat{V}_0 \bigr]^{1/3}
n^{1/3}.
\end{equation}
The bias estimator for the NPPI method is
\[
\widehat{\BIAS} = 2\bigl[\hsi^2_n(\ell_2) - \hsi^2_n(2\ell_2)\bigr]
\]
and the variance estimator is constructed using the
jackknife-after-bootstrap
(JAB) method \cite{E92,L02},
due to its computational advantages. For completeness,
we next briefly describe the details of the JAB
variance estimator.

\begin{remark}\label{rem1} Politis and Romano \cite{PR95} considered an estimator
related to $\widehat{\BIAS}$ above
for bias-correcting the Bartlett spectral estimator (e.g., at the zero
frequency, this Bartlett estimator
is asymptotically equivalent to $\hsi^2_n(\ell_2)$ and their corrected
estimator is equivalent to
$2 \hsi^2_n(2\ell_2) - \hsi^2_n(\ell_2)$). It is also important to
re-iterate that, while the NPPI
block estimator is based on general forms (cf. (\ref{eqngenexp}),
(\ref
{bv})) for the asymptotic
bias and variance of a bootstrap estimator, the HHJ block estimator
requires only the optimal block
order (cf. (\ref{eqngenblock}), (\ref{l-0})) for minimizing the
asymptotic MSE of a bootstrap estimator;
in this sense, the HHJ method requires less large-sample information
and could potentially be more general.
At the same time, as the MSE-optimal block order is typically derived
from asymptotic bias/variance
quantities, both NPPI and HHJ methods are generally intended to apply
for block selection with the
same problems, particularly under the smooth function model.
\end{remark}

%s2.3.1 #&#
\subsubsection{The jackknife-after-bootstrap variance estimator}
\label{sec231}
The JAB method was initially proposed by \cite{E92} to assess
accuracy of bootstrap estimators for independent data,
and was extended to the dependent case by \cite{L02}.
A key advantage of the JAB method is that it does \textit{not}
require a second level of resampling; the JAB method produces a
variance estimate of a block bootstrap
estimator by merely regrouping the resampled blocks
used in computing the original block bootstrap
estimator~\cite{L02}.\looseness=1

Suppose that the goal is to estimate the variance of an MBB estimator
$\hvphi_n(\ell)$ based on blocks of length $\ell$.
(For notational simplicity here, consider $\ell=\ell_1$
and $\hvphi_n(\ell) = \hsi_n^2(\ell)$.) Let $m\equiv m_n$ be
an integer such that
$
m\rightarrow\infty\mbox{ and } m/n\rightarrow0
$
as $\nti$.
Here, $m$ denotes the number of bootstrap
blocks to be deleted for the JAB.
Set $N=n-\ell+1$,
$M=N-m+1$ and for
$i=1,\ldots,M$, let $I_i=\{1,\ldots,N\} \setminus\{i,\ldots,
i+m-1\}$.
% denote the index set of all blocks of size $\ell$
% obtained by deleting the $m$ blocks
%$\{\B_i,\ldots,\B_{i+m-1}\}$.
Also, let ${\mathcal X}_{i,\ell}=(X_i,\ldots,X_{i+\ell-1})$,
$i=1,\ldots,N$
be the MBB blocks of size $\ell$.
The first step of the JAB is to define a {jackknife version}
$\hat{\vphi}_n^{(i)}\equiv\hat{\vphi}_n^{(i)}(\ell)$
of $\hat{\vphi}_n(\ell)$
for each $i\in\{1,\ldots,M\}$.
Then, the \textit{$i$th block-deleted
jackknife point value} $\hat{\vphi}_n^{(i)}$ is obtained
by resampling $ \lfloor n/\ell\rfloor$ blocks randomly, with
replacement from the reduced
collection $\{{\mathcal X}_{i,\ell} \dvt j\in I_i\}$ and then by
computing the
corresponding block bootstrap variance estimator using
the resulting resample.\looseness=1

%%%%%%%%%%%%%%%%
% Computation of the Jackknife point
%values $\hat{\vphi}_n^{(i)}$ is rather %surprisingly
% simple. Suppose that
%$\{_k{\cal B}_1^*,\ldots,_k{\cal B}_b^*\}$, $k=1,\ldots, K$
%denote the $K$ conditionally iid sets of bootstrap blocks
%drawn randomly and with replacement from the full collection
% $\{{\cal B}_1,\ldots,{\cal B}_N\}$ of overlapping blocks,
%for Monte-Carlo evaluation of the given
%block bootstrap estimator $\hat{\vphi}_n$.
%For $1\leq i\leq M$, let
%$$
%I_i^*=\Big\{k: 1\leq k\leq K, \{_k{\cal B}_1^*,
%$$
% denote the subcollection of all resampled
%block-sets $(_k{\cal B}_1^*, \ldots, _k{\cal B}_b^*)$ of size
%$b$ where none of the resampled blocks
%$_k{\cal B}_1^*,\ldots,_k{\cal B}_b^*$
%equals ${\cal B}_i, \ldots, {\cal B}_{i+m-1}$.
%Then each such $\{_k{\cal B}_1^*,\ldots,_k{\cal B}_b^*\}$
%represents a random sample from the `blocks of blocks'-deleted
%collection $\{{\cal B}_j: j\in I_i\}$ (cf. Lahiri (2002)). Thus,
%for computing $\hat{\vphi}_n^{(i)}$, we \textit{extract}
%the subcollection $\{(_k{\cal B}_1^*,\ldots,_k{\cal B}_b^*):
%k\in I_i^*\}$ from $\{(_k{\cal B}_1^*,\ldots, _k{\cal B}_b^*):
% 1\leq k\leq K\}$ and \textit{reuse} the original bootstrap
%observations to compute the bootstrap version $T_n^{*(i)}$
%of $T_n$,
%The empirical distribution of
%the replicates of
%$T_n^{*(i)}$ serves as a Monte Carlo approximation to
%$\hat{G}_{n,i}$ which, in turn, yields the $i$th jackknife point value
%$\hat{\vphi}_n^{(i)}$ via (3.5).
%
%With the jackknife point values given by (3.5),
%%%%%%%%%%%%%%%%%%%%%%%%%%%%

Then,
the JAB estimator of the variance of $\hat{\vphi}_n\equiv
\hat{\vphi}_n(\ell)$ is
given by
%e2.12 #&#
%
\begin{equation}
\widehat{\VAR}_{\JAB}(\hat{\vphi}_n)=\frac{m}{(N-m)}
\frac{1}{M}\sum_{i=1}^M \bigl(
\tilde{ \vphi}_n^{(i)}-\hat{\vphi}_n
\bigr)^2,
\end{equation}
where $\tilde{\vphi}_n^{(i)} =
m^{-1}[N\hat{\vphi}_n - (N-m)\hat{\vphi}_n^{(i)}]$
is the $i$th
\textit{block-deleted jackknife pseudo-value}
of $\hat{\vphi}_n$, $i=1,\ldots,M$.

%%%%%%%%%nmacro -for \vphi etc.

%s3 #&#
\section{Results on uniform expansion of the MSE}
\label{sec3}
%s3.1 #&#
\subsection{Assumptions}
\label{sec31}
To develop MSE and other probabilistic expansions, we
require conditions on the dependence structure of the stationary
$\mathbb{R}^d$-valued process $\{X_t\}_{t\in\mathbb{Z}}$ and the
smooth function $H$, described below. Condition \ref{coD} prescribes
differentiability assumptions on the smooth function $H$, Condition~\ref{coM}
describes mixing/moment assumptions as a function of positive
integer $r$, and Condition~\ref{coS} entails certain covariance sums are
non-zero. In particular, the sums in Condition~\ref{coS}
define the constant
$\mathcal{C}_0=[2B_0^2/V_0]^{1/3}$ in the large-sample optimal block
approximation $\ell_n^0=\mathcal{C}_0 n^{1/3}$ from (\ref{l-0}).
For
$\nu=(\nu_1,\ldots,\nu_d)\in(\mathbb{N}\cup\{0\})^d$, write
$\|\nu\|_1=\sum_{i=1}^d\nu_i$ in the following.

\renewcommand{\thecondition}{${D}$}
\begin{condition}\label{coD}
 The function
$H\dvtx\mathbb{R}^d\rightarrow\mathbb{R}$ is 3-times continuously
differentiable and $\max\{ |\partial^\nu H(x)/(\partial
x_1\cdots\partial x_d)| \dvt\|\nu\|_1=3\} \leq C(1+\|x\|^{a_0})$,
$x=(x_1,\ldots,x_d)^\prime\in\mathbb{R}^d$ for some $C>0$ and
integer $a_0\geq0$.
\end{condition}

\renewcommand{\thecondition}{${M_r}$}
\begin{condition}\label{coM}
For some $\delta>0$,
$\E\|X_1\|^{2r+\delta}<\infty$ and $\sum_{k=1}^\infty k^{2r-1}
\alpha(k)^{\delta/(2r+\delta)}<\infty$, where $\alpha(\cdot)$
denotes the strong mixing coefficient of the process
$\{X_t\}_{t\in\mathbb{Z}}$.
\end{condition}

\renewcommand{\thecondition}{${S}$}
\begin{condition}\label{coS}
$B_0 \equiv
\sum_{k=-\infty}^\infty|k|r(k) \neq0$ and $V_0\equiv
(4/3)\sigma_\infty^4 > 0$ for $\sigma_\infty^2=
\sum_{k=-\infty}^\infty r(k)$ in (\ref{sigma}), where $r(k)=
\cov(\nabla^\prime X_0, \nabla^\prime X_k)$, $k\in\mathbb{Z}$ and $
\nabla= (\partial H(\mu)/\partial x_1, \ldots,
\partial H(\mu)/\partial x_d)^\prime$ is the vector of first order
partial derivatives of $H$ at $\E X_1=\mu$.
\end{condition}

Mixing and moment assumptions as formulated in Condition \ref{coM} are standard
in investigating block resampling methods (cf. \cite{L03}, Chapter 5).
Typical expansions of the MSE of the MBB variance estimator
often require $H$ to be 2-times differentiable in the smooth function
model, whereas Condition \ref{coD}
requires slightly more in order to determine a finer expansion of this MSE.
The assumptions on the process quantities $B_0,V_0$ in Condition \ref{coS}
are mild and standard
for the block bootstrap \cite{H95,K89,L99,L07,N09,PP01,PR94}; in
particular, the assumption on $B_0$
is needed to rule out i.i.d. processes.

%s3.2 #&#
\subsection{Main results}
\label{sec32}
Recalling the MSE-approximation $f_n(\l) \equiv\l^{-2} B^2_0 +
n^{-1}\l V_0$
for the MBB variance estimator from~\eqref{mse11} (with
constants $B_0,V_0$ as in Condition \ref{coS} above),
Theorem \ref{theorem} below provides a
more refined expansion of this MSE over a collection
of block lengths, $ \mathcal{J}_n=\{\l\in\mathbb{N} \dvt
K^{-1}n^{1/3}\leq\l\leq K n^{1/3}\}$
as in \eqref{opt} (cf. Section \ref{sec21}), of optimal
order.\looseness=1

\begin{theorem}
\label{theorem} Suppose that Conditions \ref{coD}, \ref{coM} with $r=6 +2
a_0$, and Condition~\ref{coS} hold, where $a_0$ is as specified by Condition \ref{coD}.
Then, as $n \rightarrow
\infty$,
\begin{longlist}[(ii)]
\item[(i)] for $f_n(\cdot)$ defined in (\ref{mse11}),
\[
\max_{\l\in\mathcal{J}_n} \biggl\llvert\MSE_n(\l) -
\frac{2B_0
\sigma^2_\infty}{n} -f_n(\l) \biggr\rrvert= O \bigl(n^{-4/3}
\bigr).
\]
\item[(ii)] $|\l_n^{\opt} - \l_n^0|/\l_n^0=O(n^{-1/3})$,
for $\l_n^0 \equiv\operatorname{argmin}_{y>0}f_n(y) = \mathcal
{C}_0 n^{1/3}$ from (\ref{l-0}).
\end{longlist}
\end{theorem}
Theorem \ref{theorem}(i) gives a close bound $O(n^{-4/3})$ how the
MSE-approximation $f_n(\l)$
matches the curve $\MSE_n(\l) - n^{-1}2B_0 \sigma_\infty^2$
(\textit{not} quite $\MSE_n(\l)$
but both having the same minimizer), uniformly
in $\l\in\mathcal{J}_n$. For comparison, note $f_n(\l)$, $\l\in
\mathcal{J}_n$,
has exact order $O(n^{-2/3})$.
In trying to resolve $\l_n^{\opt}$, we then have a general bound on the
differences
$n^{2/3}\{\MSE_n(\l) - \MSE_n(\l_n^{\opt}) - [f_n(\l)
- f_n(\l_n^{\opt})]\} =O(n^{-2/3})$ between the two curves. One
implication, stated
in Theorem~\ref{theorem}(ii), is that
$O(n^{-1/3})$ becomes the general order on the
discrepancy between the minimizer $\l_n^{\opt} $ of $\MSE_n(\cdot)$
and the
minimizer $\l_n^0$ of $f_n(\cdot)$. \textit{Theorem \textup{\ref{theorem}}
bounds cannot
be generally improved by further
expanding $\MSE_n(\l)$} (i.e., under additional smoothness
assumptions on $H$) and, in fact in Theorem \ref{theorem}(ii), $\l
_n^{\opt}$ is necessarily an integer
while $\l_n^0$ need not be.

%%
%For example, if we assume $H$ is 4-times differentiable in
%place of condition $D$ along with additional moment conditions,
%we could strengthen the squared bias contribution to the MSE
%of the MBB variance estimator as
%[\E\hat{\sigma}_{1,m}^2(b)-\sigma_\infty^2]^2 -\left( \frac{B_0}{b}
%+ \frac{\sigma^2_\infty b}{m} \right)^2 \right| = o(m^{-4/3})\]
%and obtain a result akin to $m^{2/3}\{\MSE_m(b_m^{opt}) -
%-b_m^{\opt})/m^2) %+o(m^{-4/3})$, implying again $m^{-1/3}(b
%-b_m^{\opt})^2 = O(m^{-1/3}) $.

%s4 #&#
\section{Results on the HHJ method}
\label{sec4}
To state the main result, recall
$\hat{\ell}_{n,\HHJ}^{\opt}$ denotes the HHJ block estimator (\ref
{HHJblock}),
depending on a pilot block $\tilde{\ell}_n$ and subsample
size $m$, and that $\lhoi$ from
(\ref{HHJblock2}) denotes an
oracle-like version of $\hat{\ell}_{n,\HHJ}^{\opt}$ that requires $m$
but not
$\tilde{\ell}_n$.

\begin{theorem}
\label{theorem1} Suppose that Conditions \ref{coD}, \ref{coM} with $r=14 +4
a_0$, and Condition~\ref{coS} hold, with $a_0$ as specified by Condition \ref{coD}.
Assume that $m^{-1} + m/n\rightarrow0$ as $n\rightarrow\infty$ with
$m^{5/3}/n=O(1)$.
\begin{longlist}[(ii)]
\item[(i)] Then, as $n\rightarrow\infty$,
\[
\biggl\llvert\frac{\lhoi-\ell_n^{\opt}}{\ell_n^{\opt}} \biggr
\rrvert
=O_p \bigl( \max
\bigl\{ m^{-1/3}, m^{-1/12}(m/n)^{1/4},
(m/n)^{1/3} \bigr\} \bigr).
\]
\item[(ii)] If additionally $\tilde{\ell}_n^{-1}+ \tilde{\ell
}_n^{2}/n\rightarrow0$ and
$m/\tilde{\ell}_n^{2} + m^2/n=O(1)$, then
\[
\biggl\llvert\frac{\hat{\ell}_{n,\HHJ}^{\opt}-\ell_n^{\opt
}}{\ell
_n^{\opt}} \biggr\rrvert=O_p \biggl(\max
\biggl\{ m^{-1/3}, \frac{m^{1/3}}{\tilde{\ell}_n}, \frac{m^{1/6}}{n^{1/4}},
\frac{m^{1/3}}{(\tilde{\ell}_n n )^{1/4}}, \frac{ \tilde{\ell
}_n^{1/4}}{n^{1/4}}, \frac{m^{1/3} \tilde{\ell}_n^{1/2}}{ n^{1/2}}
\biggr\} \biggr).
\]
\end{longlist}
\end{theorem}

\begin{remark}\label{rem2}
Theorem \ref{theorem1} also holds if, on
the left-hand sides above, we replace
$(\lhoi-\ell_n^{\opt})/\ell_n^{\opt}$ and
$(\hat{\ell}_{n,\HHJ}^{\opt}-\ell_n^{\opt})/\ell_n^{\opt}$
with their subsample counterparts
$(\bhoi-b_m^{\opt})/b_m^{\opt}$
and $(\hat{b}_{m,\HHJ}^{\opt}-b_m^{\opt})/b_m^{\opt}$. This result
helps to reinforce the notion that the quality of block
estimation at the subsample level determines the performance
of HHJ method.
\end{remark}

Theorem \ref{theorem1}(i) indicates how the subsample size $m$
affects the convergence rate of the oracle-type block estimate.
It follows from Theorem \ref{theorem}(i) that, with oracle
knowledge of $\sigma_\infty^2$, the best possible (fastest)
rate of convergence for
$(\lhoi-\ell_n^{\opt})/\ell_n^{\opt}$ is $O_p(n^{-1/6})$
achieved when the subsample size $m \propto n^{1/2}$. The choice $m
\propto n^{1/2}$
balances the sizes of all three terms in
the bound from Theorem \ref{theorem1}(i). Remark \ref{rem3} below provides some
explanation of the probabilistic bounds in Theorem \ref{theorem1}(i).

In Theorem \ref{theorem1}(ii), we impose some additional
block growth conditions on the pilot block $\tilde{\ell}_n$
and subsample size
$m$ in the HHJ
method, which are mild and help to concisely express the order
of the main components contributing to the error rate.
While the combined effects of the tuning parameters are complicated
and difficult to characterize in Theorem \ref{theorem1}(ii),
a block $\tilde{\ell}_n \propto n^{1/3}$
of MSE-optimal order for the pilot
MBB variance estimator $\hat{\sigma}_n^2(\tilde{\ell}_n)$ in the HHJ
method is an intuitive starting point. And
with this choice, it follows that $m \propto n^{1/2} $
is then optimal for minimizing the convergence rate of the HHJ block
estimator, which becomes $O_p(n^{-1/6})$.
In fact, the selection $m \propto n^{1/2}, \tilde{\ell}_n\propto n^{1/3}$
is overall optimal and simultaneously balances the order
$O_p(n^{-1/6})$ of \textit{all six} error terms in
Theorem \ref{theorem}(ii). So surprisingly, the HHJ block
estimator $\hat{\ell}_{n,\HHJ}^{\opt}$ achieves the best
convergence rate that one could hope for by matching the
optimal rate of the oracle block
estimator $\lhoi$.
We summarize our findings on tuning parameters
in Corollary \ref{coro1}.

\begin{Corollary}
\label{coro1} Under the assumptions
of Theorem \ref{theorem1}, a subsample size $m \propto n^{1/2}$
and pilot block $\tilde{\ell}_n \propto n^{1/3}$ yield optimal
convergence rates
\[
\biggl\llvert\frac{\lhoi-\ell_n^{\opt}}{\ell_n^{\opt}} \biggr
\rrvert
=O_p \bigl(
n^{-1/6} \bigr),\qquad \biggl\llvert\frac{\hat{\ell}_{n,\HHJ}^{\opt
}-\ell
_n^{\opt}}{\ell_n^{\opt}} \biggr\rrvert
=O_p \bigl( n^{-1/6} \bigr)
\]
as $n \rightarrow\infty$, for the HHJ block estimator
$\hat{\ell}_{n,\HHJ}^{\opt}$ and its
oracle $\lhoi$ version.
\end{Corollary}

An interpretation of Corollary \ref{coro1} is that, at optimal
tuning parameters, random fluctuations in the HHJ block
estimator $|\hat{\ell}_{n,\HHJ}^{\opt}-\ell_n^{\opt}|$
are of the order $\sqrt{\ell_n^{\opt}}$. This behavior interestingly
resembles that of some other kernel
bandwidth estimators based empirical MSE criteria (cf.~\cite{R84}),
though $\widehat{\MSE}_m(\cdot)$
does not take its arguments from a continuum of real-values.

\begin{remark}\label{rem3} We provide a brief explanation of
the probabilistic bounds in Theorem \ref{theorem1},
and focus mainly on the behavior of oracle block estimator
$\bhoi$ from Section \ref{sec221}
at the subsample level; more rigorous details
are given in Section \ref{sec6} and the supplementary material \cite{NL12}.
Recall the block estimator
$\bhoi$ minimizes the $\widehat{\MSE}_m^{\infi}(b) $
from (\ref{mse2}),
while $b^{\opt}_m$ from (\ref{b3}) minimizes
$\MSE_m(b) \approx f_m(b)\equiv b^{-1}B_0^2 + b m^{-1}V_0$
(the subsample version
of (\ref{mse11})).
In part,
the bound $O(m^{-1/3})$ in Theorem \ref{theorem1}(i) is due
to smoothness issues with $\MSE_m(b)$ and its discrepancy from $f_m(b)$
(cf. Theorem \ref{theorem}).
The other bounds in Theorem \ref{theorem1} arise from the size of
%e4.1 #&#
%
\begin{equation}
\label{Delta}\Delta_m^{\infty}(b)= \bigl\{\widehat{
\MSE}^{\infi}_m \bigl(b_m^{\opt} \bigr) -
\E\bigl[\widehat{\MSE}^{\infi}_m \bigl(b_m^{\opt}
\bigr) \bigr] \bigr\}- \bigl\{\widehat{\MSE}^{\infi}_m(b) -\E
\bigl[\widehat{\MSE}^{\infty}_m(b) \bigr] \bigr\},
\end{equation}
$b\in\mathcal{J}_m$, where $\E[\widehat{\MSE}_m^{\infi}(b)] =
\MSE
_m(b)$; this
quantity measures the discrepancy between two differenced curves (which
should ideally match at $b=\bhoi$), where
differences
in $\MSE_m(b)$ serve to identify $b_m^{\opt}$ and similar differences
in $\widehat{\MSE}_m^{\infi}(b) $ identify $\bhoi$. It can be shown
that, for any $a_n \rightarrow0$,
\[
\max_{ b\in\mathcal{J}_m \dvt|b-b_m^{\opt}| \leq a_n m^{1/3}} a_n^{-1/2}
m^{2/3}(n/m)^{1/2} \bigl|\Delta_m^{\infty}(b)\bigr|
\]
remains stochastically bounded on shrinking neighborhoods of block
lengths around $b_m^{\opt}$,
while at the same time $\bhoi/b_m^{\opt}\stackrel{p}{\rightarrow} 1$
(i.e., $\bhoi$ is consistent for
$b_m^{\opt} \approx b_m^0=\mathcal{C}_0 m^{1/3}$); see the auxiliary result,
Theorem \ref{theorem2}, of Section \ref{sec6}. This allows other order
bounds on $(\bhoi-b_m^{\opt})/b_m^{\opt}$
to be determined by recursively ``caging'' $\bhoi$ in decreasing
neighborhoods around $b_m^{\opt}$ with high probability.
The probabilistic bounds in Theorem \ref{theorem1}(ii) are partly due
to error contributions from the MBB variance estimator
$\hat{\sigma}_n^2(\tilde{\ell}_n)$ used through $\widehat{\MSE}_m(b)$
in (\ref{mse})
to estimate $\MSE_m(b)$ at the subsample level.\looseness=1
\end{remark}

%s5 #&#
\section{Results on the NPPI method}
\label{sec5}
Next, we consider the convergence rates of the
optimal block length selector based on the NPPI
method. Recall that $r(k)= \cov(Y_1,Y_{k+1})$, $k\geq1$, where $Y_i =
\nabla'X_i
$, $i\geq1$.

\begin{theorem}
\label{thm-nppi} Suppose that Conditions \ref{coD}, \ref{coM} with $r=7 +
2a_0$, and Condition~\ref{coS} hold, with $a_0$ as specified by Condition \ref{coD}.
Assume that $\ell_2n^{-1/3}+\ell_1^{-1} +\ell_1 m^{-1} +
m/n\rightarrow0$
as $n\rightarrow\infty$.
Then, as $n\rightarrow\infty$,
%e5.1 #&#
%
\begin{eqnarray}\label{rt-np}
&& \bigl\llvert\hat{\ell}_{n,\nppi}^{\opt}-\ell_n^{\opt}
\bigr\rrvert/ \ell_n^{\opt}
\nonumber
\\[-8pt]
\\[-8pt]
\nonumber
&&\quad= O_p \bigl( [{m}/{n}]^{1/2} + [\ell_1/m] +
\ell_1^{-2} \bigr) + O_p \Biggl(
\ell_2\sum_{k=\ell_2}^{2\ell_2 -1}\bigl |r(k)\bigr|+
n^{-1/2}\ell_2^{3/2} \Biggr).
\end{eqnarray}
\end{theorem}

As the NPPI method targets the block approximation $\ell_0= [2
B_0^2/V_0]^{1/3}n^{1/3}$ (\ref{l-0}),
the first of the two terms on the right side of \eqref{rt-np}
is from the estimation of $V_0$ and the second is from
the estimation of $B_0$. For the first term, with any given choice of
$m$, the optimal choice of $\ell_1$
satisfies $\ell_1/m\pt\ell_1^{-2}$, that is, $\ell_1\pt
m^{1/3}$. For this choice of $\ell_1$, the optimal
choice of $m$ is determined by the relation $[m/n]^{1/2}\pt m^{1/3}/m$,
that is, $m\pt n^{3/7}$. Thus, the optimal rate of the first term
is $O_p(n^{-2/7})$ with $m\pt n^{3/7}$ and
$\ell_1\pt n^{1/7}$.

To determine the optimal order of the second term, first
note that the pilot block size $\ell_2$ is only required to
satisfy the constraints stated in Theorem
\ref{thm-nppi}. In particular, $\l_2$ is not required
to go to $\infty$ with the sample size.
From
%the second term on the right side of
\eqref{rt-np}, it is also evident that the
optimal choice of $\ell_2$
(to minimize the order of the second term alone)
depends on the rate of decay of
the autocovariance function $r(\cdot)$.
Since $r(k)\leq C k^{-a-1}$ for some $a\geq12$
(implied by Condition \ref{coM} with $r=7+2a_0$),
the second term can always be made
to match the optimal order of the first term,
that is, $O_p(n^{-2/7})$, by choosing
$\ell_2=O(n^{1/7})$ (note $n^{-1/2}\ell_2^{3/2}=n^{-2/7}$ in (\ref{rt-np})
when $\ell_2 \propto n^{1/7}$). However, for processes with an exponentially
decaying $r(k)$, a choice of $\ell_2 \propto\log n$ optimizes
the second term, with the attained rate of
$O_p(n^{-1/2}[\log n]^{3/2})$, while for an $m_0$-dependent sequence
$\{X_t\}$ with a fixed $m_0\geq1$, a choice of $\ell_2=m_0+1$
makes the second term $O_p(n^{-1/2})$. But, in the end, the error rate
$O_p(n^{-2/7})$ of first term dominates the second in \eqref{rt-np}.

\begin{remark}\label{rem4}
In Lahiri \textit{et al.} \cite{L07}, the NPPI
plug-in estimator
was defined with a common choice $\ell_1=\ell_2$. In this case,
under the conditions of Theorem \ref{thm-nppi}, the optimal
order of the common block size is determined by
$O_p([{m}/{n}]^{1/2} + [\ell_1/m]
+ \ell_1^{-2}+
n^{-1/2}\ell_1^{3/2})$. For a fixed $\ell_1$, the first two factors
are optimized for
\[
m\propto n^{1/3}\ell_1^{2/3}.
\]
Interestingly, this order of $m$ was also suggested by \cite{L07},
purely on the basis of some heuristic arguments. For this choice of
$m$, one may choose $\ell_1\propto n^{1/7}$ to optimize the rate
of convergence of the NPPI method, yielding the same optimal rate
$O_p(n^{-2/7})$
possible with three tuning parameters in the NPPI method. This supports
the suggestion of Lahiri \textit{et al.} \cite{L07} of a common choice
$\ell_1=\ell_2$ and, with the same number of tuning parameters, the
NPPI block selector
has a better optimal rate than $O_p(n^{-1/6})$ for the
HHJ method.
\end{remark}

We summarize our findings on the NPPI method
in Corollary \ref{coro2}.

\begin{Corollary}
\label{coro2} Under the assumptions
of Theorem \ref{thm-nppi}, a JAB block deletion size $m \propto n^{3/7}$
and tuning block lengths $\ell_1 \propto n^{1/7}$, $\ell_2=O(n^{1/7})$
as $n \rightarrow\infty$ yield an optimal convergence rate
for the NPPI block estimator
$\hat{\ell}_{n,\nppi}^{\opt}$ as
\[
\biggl\llvert\frac{\hat{\ell}_{n,\nppi}^{\opt}-\ell_n^{\opt
}}{\ell
_n^{\opt}} \biggr\rrvert=O_p
\bigl(n^{-2/7} \bigr).
\]
In particular, choosing $m \propto n^{3/7}$, $\ell_1=\ell_2 \propto
n^{1/7}$ achieves this optimal rate.
\end{Corollary}

%s6 #&#
\section{Comparison with plug-in methods and concluding remarks}
\label{sec55}
In the problem choosing an appropriate block length for implementing
block bootstraps in time series,
the HHJ and NPPI methods represent the two existing \textit{general}
block selection methods in the literature. However,
because convergence rates of these block estimators have been unknown,
our goal here was to
provide some comparison of their relative performances, considering
block estimation for MBB variance estimation
in particular. Both methods are again ``general'' in the sense that one
could consider block $\ell$ estimation for a block bootstrap
version $\hat{\varphi}_n(\ell)$ of a general functional $\varphi_n$
(e.g., bias, variance,
distribution function, quantiles, etc. as in Section \ref{sec1}) of
the sampling distribution of a time series estimator $\hat{\theta}_n$,
by replacing the MBB variance functional
$\hat{\varphi}_n(\ell) =
\hat{\sigma}_{n}^2(\ell) \equiv\ell\lfloor n/\ell\rfloor
\operatorname{Var}_*(\hat
{\theta}_n^*)$ with $\hat{\varphi}_n(\ell)$ in the mechanics
of the HHJ and NPPI methods described in Sections \ref{sec22}--\ref
{sec23}. Both methods aim to estimate
MSE-optimal block length through its
large sample approximation $\ell_n^0 = \mathcal{C}_0 n^{1/(r+2)},
\mathcal{C}_0 = [2 B_0^{2}/(r V_0)]^{1/(r+2)}$
in (\ref{eqngenblock}) and neither method requires explicit forms for
population
quantities $B_0\equiv B_0(\varphi_n), V_0\equiv V_0(\varphi_n)$
(arising, resp.,
from the bias and variance of a bootstrap functional $\hat{\varphi
}_n(\ell)$ in
(\ref{eqngenexp})) which can depend on the functional $\varphi_n$
and unknown
process parameters in a complex way. The HHJ approach estimates the constant
$\mathcal{C}_0$ in $\ell_n^0$ directly
through a subsampling technique, while the NPPI method non-parametrically
estimates both $B_0$ and $V_0$ in $\mathcal{C}_0$. Intuitively,
because the
general NPPI approach separately targets the bootstrap bias/variance $B_0,V_0$
contributions to $\mathcal{C}_0$, one might anticipate this approach
to exhibit
better convergence rates in block estimation compared to HHJ. In considering
block estimation for MBB variance estimation with time series, we have shown
that this is indeed the case. For the variance problem, NPPI achieves a
better rate $O_p(n^{-2/7})$
than the HHJ method $O_p(n^{-1/6})$ when both methods use two tuning parameters.
While considering the MBB among possible block bootstrap approaches,
the same
convergence rates and optimal tuning parameter selections should also
hold for
other block bootstraps, such as the non-overlapping block bootstrap
\cite{K89},
the circular block bootstrap \cite{PR92} and the stationary bootstrap
\cite{N09,PR94}
(though the tapered block bootstrap \cite{PP01}
requires a different treatment as the bias expansion in (\ref
{eqngenblock}) or (\ref{bv})
needs to be replaced by a smaller bias term $\ell^{-2}B_0$ in variance
estimation).
And though we have focused on variance estimation, we suspect that the
NPPI method
retains similar
large-sample superiority over the HHJ method for block selection in
other inference
problems.
%In fact, for the general functionals $\varphi_n$ in Section
%%rate $\sqrt{\ell_n^0}\propto n^{1/(2r+4)} $, with respect to $
%%seems to fit this pattern associated with kernel bandwidth estimators
%based on the variance estimation results
%in Section \ref{sec4} (e.g., a rate $O_{p}(n^{-1/6})$ for $r=1$).

As mentioned in the \hyperref[sec1]{Introduction}, in the particular
setting of block
bootstrap variance estimation, other plug-in methods for block
selection exist such as
the proposals of B\"uhlmann and K\"unsch (BK) \cite{BK99} and Politis
and White
(PW) \cite{PW04} (see also Patton, Politis and White \cite{P10}).
These use explicit expressions for the bias $B_0$ and variance
components $V_0$ of the MBB variance estimator from (\ref{bv})
appearing the approximation $\ell_n^0= [ 2 B_0^2/V_0] n^{1/3}$ of the
MSE-optimal block length $\ell_n^{\opt}$ (\ref{opt}), given in this
case by
%e6.1 #&#
%
\begin{equation}
\label{eqnBV} B_0 = \sum_{k=-\infty}^\infty
|k| r(k), \qquad V_0 = \frac{4}{3} \Biggl( \sum
_{k=-\infty}^\infty r(k) \Biggr)^4
\end{equation}
for $r(k)= \cov( \nabla'X_1, \nabla'X_{1+k})$, $k\geq1$; see Condition
\ref{coS}, Section \ref{sec31}.
The BK and PW approaches estimate the covariance sums (\ref{eqnBV})
with spectral lag window estimators which are then plugged into the
approximation $\ell_n^0$
(\ref{l-0}) to estimate $\ell_n^{\opt}$.
The BK method is based on an iterative plug-in algorithm from
B\"uhlmann \cite{B96}
for estimating the optimal bandwidth for lag window estimators
of the spectral density at zero, which has equivalences to block length
selection
for the MBB variance estimator.
If $\hat{\ell}_{n,\BK}^{\opt}$ denotes the resulting block estimator,
results in
B\"uhlmann and K\"unsch \cite{BK99} show that
\[
\frac{\hat{\ell}_{n,\BK}^{\opt} - \ell_n^{\opt}}{\ell_n^{\opt
}} = O_p
\bigl(n^{-2/7} \bigr)
\]
for MBB variance estimation of smooth function model statistics. As
mentioned earlier,
interestingly, the NPPI block selection method
obtains the exact \textit{same} optimal rate of convergence without
using the structural
knowledge in (\ref{eqnBV}). The plug-in estimator $\hat{\ell
}_{n,\PW
}^{\opt}$ of
Politis and White \cite{PW04} is formulated by using a ``flat-top''
lag-window $ \lambda(t) =
\mathbb{I}(t\in[0,1/2])+ 2(1-|t|)\mathbb{I}(t \in(1/2,1])$, $t \in
[0,1]$, where
$\mathbb{I}(\cdot)$ denotes the indicator function; see \cite{PR95}.
Their method
was originally studied for block bootstrap variance estimation of time series
sample means. Here we describe an extension of the methodology
%can be modified
for smooth function model statistics
$\hat{\theta}_n=H(\bar{X}_n)$. The corresponding two
unknown covariance sums (\ref{eqnBV}) in $\ell_n^0$ are
estimated, respectively, with
%e6.2 #&#
%
\begin{equation}
\label{eqnPW}\sum_{k=-2M}^{2M} \lambda
\bigl\{k/(2M) \bigr\} |k| \hat{\nabla} \hat{r}(k) \hat{\nabla
}^\prime,\qquad
\sum_{k=-2M}^{2M} \lambda\bigl\{k/(2M) \bigr
\}\hat{ \nabla} \hat{r}(k) \hat{\nabla}^\prime,
\end{equation}
where
$\hat{r}(k)=n^{-1}\sum_{i=1}^{n-|k|}
(X_i-\bar{X}_n)(X_{i+|k|}-\bar{X}_n)^\prime$, $\hat{\nabla
}=\partial
H(\bar{X}_n)/\partial x$, and $M$ is a positive integer bandwidth. In
which case, we may
state a result on the convergence rate of the generalized PW block estimator.
\begin{theorem}
\label{theorem3} Under the assumptions
of Theorem \ref{theorem}, if $M\propto n^\tau$ as $n\rightarrow
\infty$
for some $10^{-1}\leq\tau\leq3^{-1}$,
Politis--White block estimator
$\hat{\ell}_{n,\PW}^{\opt}$ satisfies
\[
\biggl\llvert\frac{\hat{\ell}_{n,\PW}^{\opt}-\ell_n^{\opt
}}{\ell
_n^{\opt}} \biggr\rrvert=O_p
\bigl(n^{-1/3} \bigr).
\]
This also holds for other rules for selecting $M$ under Theorem 3.3
conditions of \cite{PW04}.\vspace*{-1pt}
\end{theorem}

\begin{remark}\label{rem5}
It should be noted that the rate in Theorem \ref{theorem3}
differs from
results in Politis and White~\cite{PW04} who considered a different
problem in
block estimation. Namely,\vspace*{1pt} they considered convergence rates between
$\ell_n^0= [ 2 B_0^2/V_0] n^{1/3}$ and
its plug-in counterpart $\hat{\ell}_{n,\PW}^{\opt}= [ 2 \hat
{B}_0^2/\hat
{V}_0] n^{1/3}$,
where $\ell_n^0$ again represents the large-sample
approximation (\ref{l-0}) of the MSE-optimal block length $\ell
_n^{\opt
}$ from (\ref{opt}).
They showed that, depending on
the underlying process dependence (cf. their Theorem 3.3), the
bandwidth $M$ can be
adaptively chosen so that $|\hat{\ell}_{n,\PW}^{\opt}
- \ell_n^0|/\ell_n^0$ may exhibit a convergence rate as high as
$O_p(n^{-1/2})$;
see also Politis \cite{P03} for a related
discussion of rate adaptivity and empirical rules for selecting $M$.
In these cases, there is still a bound $O(n^{-1/3})$
on the relative closeness of $\ell_n^0$ and $\ell_n^{\opt}$ from
Theorem \ref{theorem}.
Additionally, while the PW and NPPI methods both involve plug-in
estimation, the NPPI
approach does not require or use an explicit form
for $B_0,V_0$ in the variance problem (\ref{eqnBV}), and the
discussion of Section \ref{sec5}
indicates that this method can adaptively estimate $B_0$ (with
similar rates as high as $O_p(n^{-1/2})$) but does not adaptively
estimate $V_0$. That is,
the JAB (i.e., block jackknife) variance estimator for $V_0$ is not
rate adaptive in the NPPI method,
but the PW flat-top kernel approach is. These differences explain the
superior performance of the PW
method compared to NPPI for block estimation in the variance estimation
problem.\vspace*{-1pt}
\end{remark}
%t1 #&#
%
\begin{table}[b]
\tabcolsep=0pt
\caption{Optimal convergence rate $|\hat{\ell}_n^{\opt} -\ell
_n^{\opt}
|/ \ell_n^{\opt}$ for block estimators
$\hat{\ell}_n^{\opt}$ of the MSE-optimal block length $\ell_n^{\opt}$
(\protect\ref{opt})
for MBB variance estimation, based on the approximation $\ell_n^0$
(\protect\ref{l-0})}\label{tab1}
\begin{tabular*}{\textwidth}{@{\extracolsep{\fill}}llllll@{}}
\hline
& \multicolumn{5}{l@{}}{Methods}\\[-6pt]
& \multicolumn{5}{l@{}}{\hrulefill}\\
&\multicolumn{2}{l}{General} & \multicolumn{2}{l}{Form (\ref
{eqnBV})-based plug-in} &
\\[-6pt]
&\multicolumn{2}{l}{\hrulefill} & \multicolumn{2}{l}{\hrulefill} &
\\
& \multicolumn{1}{l}{HHJ} & \multicolumn{1}{l}{NPPI} & \multicolumn
{1}{l}{B\"uhlmann--K\"unsch (BK)}& \multicolumn{1}{l}{Politis--White (PW)}&
\multicolumn{1}{l@{}}{Best possible}\\%[-6pt]
%& \multicolumn{1}{l}{\hrulefill} & \multicolumn{1}{l}{\hrulefill} &
\hline
Rate & $O_p(n^{-1/6})$ & $O_p(n^{-2/7})$ & $O_p(n^{-2/7})$ &
$O_p(n^{-1/3})$ &
$O_p(n^{-1/3})$\\
\hline
\end{tabular*}
\end{table}

Table \ref{tab1} provides a final summary of the convergence rates of both
general and (\ref{eqnBV})-based
plug-in methods for block selection with the MBB variance estimator.
The PW plug-in estimator
attains the highest convergence rate $O_p(n^{-1/3})$ possible under
Theorem \ref{theorem} for
any estimator of MSE-optimal block length $\ell_n^{\opt}$ which is
based on its asymptotic
approximation $\ell_n^0$ (\ref{l-0}). That is, the
plug-in method of Politis and White \cite{PW04} has the \textit{best}
large-sample properties
of any existing method for block selection in the variance estimation
problem with mean-like\vadjust{\goodbreak}
or smooth function model statistics. Of course, this advantage comes at
the price in that the PW method is
designed for variance estimation (i.e., the forms (\ref{eqnBV}) in
this problem) and is therefore ``non-general''
or not directly usable for block selection in other block bootstrap
applications.
In particular, for other inference problems (e.g., distribution or
quantile estimation),
the forms of $B_0,V_0$
in the large-sample block formulas (\ref{eqngenblock}) can become
complicated, depending
additionally sums of higher order process cumulants in a more
complex fashion than the variance estimation problem. In these cases,
where appropriate
block selections
for the block bootstrap are still needed, the general HHJ and NPPI
block estimation methods
have their greatest appeal, and the convergence rate results in
variance estimation suggest
that the NPPI may have better performance than HHJ more generally.

%s7 #&#
\section{Additional results and proofs}
\label{sec6}

Theorem \ref{lem1} below gives
a bias and variance decomposition for the MBB variance
estimator $\hat{\sigma}^2_{m}(b)$, uniformly in $b\in\mathcal{J}_m$,
which is used to establish Theorem \ref{theorem}.
The proof of Theorem \ref{lem1} appears in the supplementary material
\cite{NL12}.

\begin{theorem}
\label{lem1} Under the assumptions of Theorem \ref{theorem}, as $m
\rightarrow
\infty$,
\begin{eqnarray*}
\mathrm{(i)} &&\qquad  \max_{b\in\mathcal{J}_m} \biggl
\llvert
\bigl[\E\hat{\sigma}^2_{m}(b) -\sigma^2_\infty
\bigr] + \biggl( \frac{B_0}{b} + \frac{b}{m}\sigma_\infty^2
\biggr) \biggr\rrvert= O \bigl(m^{-1} \bigr),
\\
\mathrm{(ii)}&&\qquad \max_{b\in\mathcal{J}_m} \biggl\llvert
\operatorname{Var}\bigl[ \hat{
\sigma}^2_{m}(b) \bigr] - V_0
\frac{b}{m} \biggr\rrvert= O \bigl(m^{-4/3} \bigr).
\end{eqnarray*}
\end{theorem}

\begin{pf*}{Proof of Theorem \ref{theorem}} Part (i) follows directly from
Theorem \ref{lem1}. Part (ii) follows by
expanding $0\leq n^{2/3}[\MSE_n(\lfloor\ell_n^0\rfloor) - \MSE
_n(\ell
_n^{\opt})]$ with Theorem \ref{theorem}(i),
implying $0 \leq n^{2/3}[f_n(\ell_n^{\opt})-f_n(\ell_n^0)]\leq C
n^{-2/3}$. Then, a second order Taylor expansion of
$f_n(\cdot)$ around $\ell_n^0=\mathcal{C}_0n^{1/3}$ gives the result
(as $d f_n( \ell_n^0)/dy =0$).
\end{pf*}

Theorem \ref{theorem2} next establishes the consistency of the HHJ
block estimator (and its oracle-version)
at both sample and subsample levels, and provides tightness results for
developing rates for the HHJ block estimator (cf. Theorem \ref
{theorem1}); its
proof is given in the supplementary material \cite{NL12}.
To state the result, recall the difference $\Delta_m^{\infi}(b)$,
$b\in
\mathcal{J}_m$, between empirical and true MSE curves from
(\ref{Delta}) and define $\Delta_m(b)$ by replacing $\widehat{\MSE
}^{\infi}_m(\cdot)$ from (\ref{mse2}) with
$\widehat{\MSE}_m(\cdot)$ from (\ref{mse}) in (\ref{Delta}) (i.e., HHJ
method uses $\widehat{\MSE}_m$). Then,
%e7.1 #&#
%
\begin{equation}
\label{Delta2} \Delta_m(b) = \Delta_m^{\infi}(b)
+\Omega_{1,m}(b) + \Omega_{2,m}(b),\qquad \Omega_{1,m}(b)
\equiv\Omega_{3,m}(b)-\E\bigl[\Omega_{3,m}(b) \bigr]
\end{equation}
holds for $b \in\mathcal{J}_m$, where $\Omega_{2,m}(b)\equiv2[\hat
{\sigma}^2_{n}(\tilde{\ell}_n) -
\E\hat{\sigma}^2_{n}(\tilde{\ell}_n)] \E[\hat{\sigma
}^2_{1,m}(b_m^{\opt
}) - \hat{\sigma}^2_{1,m}(b)] $ and
\[
\Omega_{3,m}(b)\equiv2\frac{\hat{\sigma}^2_{n}(\tilde{\ell}_n)
-\sigma
_\infty^2}{n-m+1}\sum
_{i=1}^{n-m+1} \bigl\{ \bigl[\hat{\sigma}^2_{i,m}(b)
- \hat{\sigma}^2_{i,m} \bigl(b_m^{\opt}
\bigr) \bigr] -\E\bigl[\hat{\sigma}^2_{i,m}(b) - \hat{
\sigma}^2_{i,m} \bigl(b_m^{\opt} \bigr)
\bigr] \bigr\}.
\]
Given any $C>0$, define a
block set $\mathcal{J}_m^{\opt}(C) = \{b\in\mathcal{J}_m\dvt
|b_m^{\opt
}-b|\leq C
m^{1/3}\}$.\vadjust{\goodbreak}

\begin{theorem}
\label{theorem2} Suppose that Conditions \ref{coD}, \ref{coM} with $r=14 +4
a_0$, and Condition~\ref{coS} hold, with $a_0$ as specified by Condition \ref{coD}.
Assume that $m^{-1}+m/n\rightarrow0$ with $m^{5/3}/n=O(1)$ as
$n\rightarrow\infty$ and that
$\tilde{\ell}_n$ in the HJJ method satisfies $ \tilde{\ell}_n^{-1} +
\tilde{\ell}^2_n/n \rightarrow0$ and $ m(\tilde{\ell}_n^{-2} +
n^{-1}\tilde{\ell}_n) =O(1)$.
Let
$\Lambda_m(b)$, $b\in\mathcal{J}_m$ denote either $\Delta_m^{\infi}(b)$
or $\Delta_m^{\infi}(b) + \Omega_{1,m}(b)$.
Then,
\begin{longlist}[(iii)]
\item[(i)] there exists an integer $N_0 \geq1$ and constant
$A>0$ such that
\begin{eqnarray*}
P \biggl( a_n^{-1/2}m^{2/3} \biggl(
\frac{n}{m} \biggr)^{1/2}\max_{b\in\mathcal{J}_m^{\opt}(a_n)}
\bigl
\llvert\Lambda_m(b) \bigr\rrvert>\lambda\biggr) &\leq&
\frac{A}{\lambda},
\\
P \biggl( a_n^{-1} m^{2/3} \biggl(
\frac{m^{1/3}}{\tilde{\ell}_n} \frac{n}{m} \biggr)^{1/2}\max
_{b\in\mathcal{J}_m^{\opt}(a_n)} \bigl\llvert\Omega_{2,m}(b)
\bigr
\rrvert>
\lambda\biggr) &\leq&\frac{A}{\lambda},
\end{eqnarray*}
holds for any $\lambda>0$, any $n \geq N_0$ and any positive
$ a_n >0$.

\item[(ii)] $\hat{b}_{m,\HHJ}^{\opt}/ b_m^{\opt}
\stackrel{p}{\rightarrow}1$ and $\bhoi/ b_m^{\opt}
\stackrel{p}{\rightarrow}1$ as $n \rightarrow\infty$.

\item[(iii)] $\hat{\ell}_{n,\HHJ}^{\opt}/ \ell_n^{\opt}
\stackrel{p}{\rightarrow}1$ and $\lhoi/ \ell_n^{\opt}
\stackrel{p}{\rightarrow}1$ as $n \rightarrow\infty$.
\end{longlist}
\end{theorem}

\begin{pf*}{Proof of Theorem \ref{theorem1}} We establish Theorem \ref
{theorem1}(i) here and defer
the proof of Theorem \ref{theorem1}(ii) to the supplementary material
\cite{NL12}.
For the minimizer $b_m^{\opt}$ of $\MSE_m(\cdot)$ from (\ref{b3})
and the minimizer $b_m^0 =
\mathcal{C}_0 m^{1/3}$ of $f_m(y)$, $y>0$
(i.e., subsample version of (\ref{l-0}) solving $d [f_m(b_m^0)]/dy=0$),
Theorem \ref{theorem}(ii) gives $m^{-1/3}|b^0_m-b_m^{\opt}|=O(m^{-1/3})$,
$m^{2/3}|f_m(b^0_m)-f_m(b_m^{\opt})|=O(m^{-2/3})$ so that
%e7.2 #&#
%
\begin{equation}
\label{exmse} 0 \leq m^{2/3} \bigl[f_m\bigl(
\bhoi\bigr)-f_m \bigl(b_m^0 \bigr) \bigr] \leq C
m^{-2/3} +m^{2/3} \bigl[\MSE_m\bigl(\bhoi\bigr) -
\MSE_m \bigl(b_m^{\opt} \bigr) \bigr],
\end{equation}
by Theorem \ref{theorem}(i), for a constant $C>0$ independent of $m$.
Applying a Taylor
expansion of $f_m(\bhoi)$ around $b_m^0$ and Theorem \ref{theorem}(ii),
there exists a constant $C_0>0$ for which
\[
m^{-2/3} \bigl( \bhoi-b_m^{\opt}
\bigr)^2 \leq C_0\max\bigl\{m^{-2/3},
m^{2/3} \bigl[\MSE_m\bigl(\bhoi\bigr) - \MSE_m
\bigl(b_m^{\opt} \bigr) \bigr] \bigr\},
\]
whenever $|\bhoi/b_m^0 -1|<1/2$. Also, by definition
we have
%e7.3 #&#
%
\begin{equation}
\label{expand3} 0 \leq m^{2/3} \bigl[\MSE_m\bigl(\bhoi\bigr) -
\MSE_m \bigl(b_m^{\opt} \bigr) \bigr] \leq
m^{2/3} \Delta^{\infi}_m\bigl(\bhoi\bigr)
\end{equation}
for $\Delta^{\infi}_m(\cdot)$ defined in (\ref{Delta}) where $\E
[\widehat{\MSE}^{\infi}_m(b)]=\MSE_{m}(b)$, $b\in\mathcal{J}_m$
so that
%e7.4 #&#
%
\begin{equation}
\label{step1} m^{-2/3} \bigl( \bhoi-b_m^{\opt}
\bigr)^2 \leq C_1\max\bigl\{m^{-2/3},
m^{2/3}\Delta^{\infi}_m\bigl(\bhoi\bigr) \bigr\},
\end{equation}
must hold for any $C_1 > C_0$ whenever $|\bhoi/b_m^0
-1|<1/2$; since this last event has arbitrarily large probability by
Theorem \ref{theorem2}(ii), we will always assume (\ref{step1}) to
hold without
loss of generality along with $\mathcal{J}_m^{\opt}(C_0)=\mathcal
{J}_m$, defining a block set
$\mathcal{J}_m^{\opt}(C) = \{b\in\mathcal{J}_m\dvt|b_m^{\opt
}-b|\leq C
m^{1/3}\}$ for any $C>0$.

We next formulate a series of recursive events to coerce
$\bhoi$ into shrinking neighborhoods around $b_{m}^{\opt}$ with high
probability.
Fix $C_1>C_0$ and define $a_{0,n}=C_1$, $L \equiv\lceil\log\log n
\rceil>1$, and
\[
a_{i,n}^2\equiv C_1^2\max\biggl
\{m^{-2/3}, 2^{\sum_{k=0}^{i-1} [(L-i+1)+k]4^{-k}} \biggl(\frac{m}{n}
\biggr)^{2^{-1}\sum_{k=0}^{i-1} 4^{-k}} \biggr\},\qquad i\geq1.
\]
%
% and $m$ satisfies $(\log
%n)^2(m/n)\rightarrow0$ as $n\rightarrow\infty$ by assumption (which
%ensures
%$\max_{1\leq i \leq L} a_{i,n}^2\rightarrow0$).
Define an integer
$
J = \min\{i=1,\ldots, L +1\dvt a_{i,n}^2 = C_1^2 m^{-2/3}
\}$
and let $J=L +1$ if this integer set is empty.
For $i=0,\ldots,J-1$, let $A_i$ be the event
\[
\max_{b \in\mathcal{J}_m^{\opt}(a_{i,n})} m^{2/3}\bigl|\Delta^{\infi}_m(b)\bigr|
\leq a_{i,n}^{1/2} (m/n)^{1/2} \lambda_i,\qquad
\lambda_i \equiv C^{1/2}_1 2^{L-i}
\]
and let $B_{i}$, $i \geq1$, be the event
\[
m^{-2/3} \bigl( \bhoi-b_m^{\opt}
\bigr)^2 \leq a_{i,n}^2.
\]
Since $\mathcal{J}^{\opt}_m(a_{0,n})=\mathcal{J}_m$, event $A_0$ implies
$B_1$ by (\ref{step1}). Also, for $J>1$, if $A_i \cap B_i$ holds for
some $i=1,\ldots,J-1$,
then so must $B_{i+1}$ by (\ref{step1}), which in turn implies $\bhoi$
in the block neighborhood
$\mathcal{J}_m^{\opt}(a_{i+1,n})$ for event $A_{i+1}$.
Suppose now that $A_{J}\cap B_{J}$ holds for an event $A_{J}$
defined as
\[
\max_{b \in\mathcal{J}_m^{\opt}(a_{J,n})} m^{2/3}\bigl|\Delta^{\infi}_m(b)\bigr|
\leq a_{n,J}^{1/2} (m/n)^{1/2} C_1^{1/2};
\]
the complement $(A_{J}\cap B_{J})^c$ has probability bounded by
$\sum_{i=0}^{J} P(A^c_i) \leq A C_1^{-1/2}(1 +  \sum_{k=0}^L 2^{-k})
\leq3 A C_1^{-1/2}$ by Theorem \ref{theorem2}(i),
which can be made arbitrarily small by large $C_1$. When $A_{J}\cap
B_{J}$ holds, then by construction (\ref{step1}) further implies
that either
\[
m^{-2/3} \bigl( \bhoi-b_m^{\opt}
\bigr)^2 \leq C_1^2\max\bigl
\{m^{-2/3}, m^{-1/6}(m/n)^{1/2} \bigr\}
\]
if $a_{J,n}^2=C_1^2 m^{-2/3}$, or the remaining possibility is
$a_{J,n}\neq C_1^2 m^{-2/3}$ in event $A_J$ and $J=L +1$ so that
\begin{eqnarray*}
m^{-2/3} \bigl( \bhoi-b_m^{\opt}
\bigr)^2 & \leq& C_1^2\max\biggl
\{m^{-2/3}, 2^{\sum_{k=0}^{L} k 4^{-k}} \biggl(\frac{m}{n}
\biggr)^{2^{-1}\sum_{k=0}^{L} 4^{-k}} \biggr\}
\\
&\leq& 2^{3/2}C_1^2\max\biggl
\{m^{-2/3}, \biggl(\frac{m}{n} \biggr)^{2/3} \biggr\}
\end{eqnarray*}
using that $ (m/n)^{2^{-1}\sum_{k=0}^{L}
4^{-k}-2/3}\leq2$ and $2^{\sum_{k=1}^{L} k 4^{-k}}\leq\sqrt{2}$
for $L=\lceil\log\log n \rceil$. Hence,
\[
m^{-2/3} \bigl( \bhoi-b_m^{\opt}
\bigr)^2 \leq4 C_1^2 \max\bigl\{
m^{-2/3}, m^{-1/6}(m/n)^{1/2}, (m/n)^{2/3}
\bigr\}
\]
holds with arbitrarily high probability (large $C_1$).
Because $|b_m^{\opt} - b_m^0|m^{-1/3}=O(m^{-1/3})$ and $|\ell_n^{\opt}
- \ell_n^0|n^{-1/3}=O(n^{-1/3})$
by Theorem \ref{theorem}(ii), where $\ell_n^0 \equiv\mathcal{C}_0
n^{1/3} = b_m^0 (n/m)^{1/3}$ and
$\lhoi= \bhoi(n/m)^{1/3}$ is formed by
rescaling (\ref{HHJblock2}), Theorem \ref{theorem1}(i) follows.
\end{pf*}

\begin{pf*}{Proof of Theorem \ref{thm-nppi}}
We sketch the proof, providing more technical detail in the
supplementary material \cite{NL12}.
Considering $\hvo$ and letting $p=n/\ell_1$, it can be shown that
%e7.5 #&#
%
\begin{equation}
\bigl| \widehat{\VAR} -\operatorname{Var}\bigl(\hat{\sigma}_n^2(
\ell_1) \bigr)\bigr | = O_p \bigl(p^{-1}
\bigl[[{m}/{n}]^{1/2} % + {m}/{n}
+ \ell/m \bigr] \bigr). \label{vhat}
\end{equation}
Next using the arguments in the proof of Theorem \ref{lem1}(ii),
one can show that
%e7.6 #&#
%
\begin{equation}
\label{vhat-b} \operatorname{Var}\bigl(\hat{\sigma}_n^2(
\ell_1) \bigr) = \frac{V_0}{p}+ O \bigl(\ell_1^{-1}n^{-1}
\bigr).
\end{equation}
Hence, by \eqref{vhat} and \eqref{vhat-b}, it follows that
\begin{eqnarray*}
|\hat{V}_0 - V_0| &=& \bigl|n\ell_1^{-1}
\widehat{\operatorname{VAR}} - V_0 \bigr|
\\
&\leq& \bigl|n\ell_1^{-1} \bigl(\widehat{\operatorname{VAR}} -
\operatorname{Var}\bigl(\hat{\sigma}_n^2(\ell_1) \bigr)
\bigr)\bigr | + \bigl|n\ell_1^{-1}\operatorname{Var}\bigl(\hat{
\sigma}_n^2(\ell_1) \bigr) - V_0
\bigr|
\\
%&=& O\Big(\frac{n}{\ell} \Big)
% + O_p([{m}{n}]\cdot p^{-1})
% + O_p([\ell/m]\cdot p^{-1}])\Big] +O(n^{-1/3}\ell^{-1}) \\
&=& O_p
\bigl([{m}/{n}]^{1/2} + [\ell_1/m] +\ell_1^{-2}
\bigr).
\end{eqnarray*}
Next consider $\hbo$.
%By Taylor's expansion around
%$\hmun\equiv\hmun(\ell_1) =\E_*\bXn$ (cf. \eqref{LQC}), we have
%n_1b\var_*(L_{n}^*) + n_1b\var_*(Q_{n}^*) +
%n_1b\var_*(C_{n}^*) \\
%2n_1b\cov_*(Q_n^*,C_n^*)
%where $n_1= \lfloor n/\ell_1\rfloor\ell_1$ and where $L_n^*,Q_n^*,
%C_n^*$
%are defined as in \eqref{LQC}
% but using the MBB with block size $\ell_1$ and sample size $n$.
Using arguments in the proof of Theorem \ref{lem1}, one can show that
\[
\E\hsi_n^2(k) = k \E[\bar{Y}_k]^2
+O \bigl(n^{-1}k \bigr),\qquad \operatorname{Var}\bigl(\hsi_n^2(k)
\bigr) = O \bigl(n^{-1}k \bigr)
\]
for $k= \ell_2, 2\ell_2$, where $\bar{Y}_k = k^{-1}
\sum_{i=1}^k Y_i$ and $Y_i = \nabla'X_i$, $i\geq1$. Hence, it
follows that
\begin{eqnarray*}
\hbo&=& 2\ell_2 \bigl(\ell_2\E[\bar{Y}_{\ell_2}]^2-2
\ell_2\E[\bar{Y}_{2\ell_2}]^2 \bigr) +
O_p\bigl(n^{-1/2}\ell_2^{3/2}\bigr)
\\
&=& B_0+ O \Biggl(\sum_{k=\ell_2}^{2\ell_2 -1}
k \bigl|r(k)\bigr| + \ell_2 \sum_{k=\ell_2}^{2\ell_2 -1}
\bigl|r(k)\bigr| \Biggr) + O_p\bigl(n^{-1/2}\ell_2^{3/2}
\bigr).
\end{eqnarray*}
Combining the bounds on $\hvo$ and $\hbo$, the theorem
follows.
\end{pf*}

\begin{pf*}{Proof of Theorem \ref{theorem3}} From the assumed conditions,
$\hat{\nabla}\equiv\partial
H(\bar{X}_n)/\partial x =\nabla+ O_p(n^{-1/2})$ holds for $\nabla
\equiv\partial H(\mu)/\partial x$.
By (\ref{eqnBV})--(\ref{eqnPW}), the PW
block estimator for MBB variance estimation can be written as
$n^{-1/3}\hat{\ell}_{n,\PW}^{\opt} =
n^{-1/3}\tilde{\ell}_{n,\PW} + O_p(n^{-1/2})$ for
\[
\tilde{\ell}_{n,\PW}= (3/2)^{1/3} \Biggl( \sum
_{k=-2M}^{2M} \lambda\bigl\{k/(2M) \bigr\} |k|
\hat{r}_Y(k) \Biggr)^{2/3} \Biggl( \sum
_{k=-2M}^{2M} \lambda\bigl\{k/(2M) \bigr\}
\hat{r}_Y(k) \Biggr)^{-4/3} n^{1/3}
\]
with $\hat{r}_Y(k)=n^{-1}\sum_{i=1}^{n-|k|}
(Y_i-\bar{Y}_n)^2$ with $Y_i = \nabla^\prime X_i$, $i \geq1$. By
Theorem 3.3(i) of \cite{PW04}
with $M\propto n^\tau$ for some $10^{-1}\leq\tau\leq3^{-1}$ and the
assumed mixing conditions here,
$|\tilde{\ell}_{n,\PW}^{\opt} - \ell_n^0|/\ell_0
=O_{p}(n^{-(1-\tau)/2
}) =O_p(n^{-1/3})$ follows
for $\ell_n^0= [2 B_0^2/V_0]^{1/3} n^{1/3}$ with $B_0,V_0$ in\vadjust{\goodbreak} (\ref
{eqnBV}). Hence, by
Theorem\ref{theorem}(ii) and $\ell_n^0/\ell_n^{\opt}\rightarrow1$,
\begin{eqnarray*}
\frac{|\hat{\ell}_{n,\PW}^{\opt} - \ell_n^{\opt} |}{\ell
_n^{\opt}} &\leq
& \frac{|\hat{\ell}_{n,\PW}^{\opt}-\tilde{\ell}_{n,\PW} |}{\ell
_n^{\opt
}} + \frac{|\tilde{\ell}_{n,\PW} -\ell_n^0 |}{\ell_n^{\opt}}+
\frac{| \ell_n^0-\ell_n^{\opt} |}{\ell_n^{\opt}}
\\
&=&O_p \bigl(n^{-1/2} \bigr) + O_p
\bigl(n^{-1/3} \bigr)+O \bigl(n^{-1/3} \bigr)=O_p
\bigl(n^{-1/3} \bigr).
\end{eqnarray*}
\upqed\end{pf*}

\section*{Acknowledgments}
The authors are grateful to two referees for thoughtful comments and suggestions
which clarified and improved the manuscript, and also led to the
results in Section \ref{sec55}.
Soumendra N. Lahiri
supported in part by NSF Grant no. DMS-10-07703 and
NSA Grant no. H98230-11-1-0130.
Daniel J. Nordman
supported in part by NSF Grant no. DMS-09-06588.

\begin{supplement}%[id=suppA]
\stitle{Proofs of main results for empirical
block length selectors\\}
\slink[doi]{10.3150/13-BEJ511SUPP} %[doi,text={...}] - jei reikia
%suskaldyti doi
\sdatatype{.pdf}
\sfilename{BEJ511\_supp.pdf}
\sdescription{A supplement \cite{NL12} provides more detailed proofs of
the main
results (Theorems \ref{theorem1}--\ref{thm-nppi})
about the convergence rates for the HHJ/NPPI block selection methods
from Sections~\ref{sec4}--\ref{sec5},
as well as proofs for the auxiliary results
(Theorems~\ref{lem1}--\ref{theorem2}) of Section \ref{sec6}.}
\end{supplement}

% imsref loaded by akundreckaite, 2013-03-25 17:30:39
%

\printhistory

\end{document}